On Hypergeometric 3F2(1)  - A Review

Michael Milgram[*], Consulting Physicist, Geometrics Unlimited, Ltd.,
Box 1484, Deep River, Ont., Canada.  K0J 1P0

**Abstract:** By systematically applying ten well-known and inequivalent two-part relations between hypergeometric sums $_3F_2(\ldots|1)$ to the published database of all such sums, 62 new sums are obtained.  The existing literature is summarized, and many purportedly novel results extracted from that literature are shown to be special cases of these new sums.  The general problem of finding elements contiguous to Watson's, Dixon's and Whipple's theorems is reduced to a simple algorithm suitable for machine computation.  Several errors in the literature are corrected or noted. The present paper both summarizes and extends a previous work on this subject.

## 1. Introduction

The evaluation of the hypergeometric sum $_3F_2(\ldots|1)$ is of ongoing interest, since it appears ubiquitously in many physics and statistics problems.  Particularly tantalizing is the fact that Gauss' theorem gives a simple result for $_2F_1(\ldots|1)$, and Bühring[1] has classified general limiting cases of $_3F_2(\ldots|1)$ with distinct parametric excess, although Wimp[2] has demonstrated that no representation consisting solely of gamma functions exists for the general case.  With reference to the WZ *certification* algorithm[3] which requires that the candidate sum contain a positive integer parameter *"n"*, one might expect that computer algebra systems would by now have become a repository of existing information about such sums, but attempts to extract evaluations for particular values of $_3F_2(\ldots|1)$ rarely generate a useful outcome (see Section 3), in spite of the availability of Petkovšek's implementation[4] of the WZ approach.  In the case of non-terminating series, Koornwinder has shown[5] that the WZ method can be used to obtain 3-part transformation identities, but again it is rare that a computer algebra program will generate a closed-form result for a sum of interest. Thus, in practice, the practitioner is usually reduced to scanning published tables or the literature in the hope that a particular problem at hand has a known closed-form[6] evaluation.

At the same time, it is commonly unrecognized that the potential universe of closed-form results is much larger than those listed in the standard tables[7], because of the Thomae relations that couple two $_3F_2(\ldots|1)$ with transformed parameter sets.  Since there are 10 inequivalent forms of the 120 Thomae transformations, (see Appendix A), the potential universe of "knowable" results is nine times larger than that given in the usual tables (one of the 10 relations is the identity).  With the advent of computer algebra systems, it has now become possible to explore this universe for potentially new, closed-form results, keeping in mind the dictum[3] that such a database can never be complete.

This paper describes such a study and places it in the context of known results. In Section 2 the algorithmic details are given. Section 3 gives general background on specific closed-form results that were examined and records a number of errors in, and corrections for, items that appear in the literature (and tables).  In Section 4 observations are made that lead to two sums that were previously conjectured. Appendix A lists the ten inequivalent Thomae relations used and Appendix B lists some selected results from a previous work[8]. Appendix C gives some particular results that are useful for obtaining closed forms contiguous to Watson's, Whipple's and Dixon's theorems (see Section 3), as well as contiguity relations that permit the evaluation of a new set of sums identified in Section 5. In a preliminary version of this work[8], a large number of useful relations were identified that transform a candidate infinite series $_3F_2(\ldots|1)$ involving a positive integer parameter *"n"* into a finite series $_3F_2(-n,\ldots|1)$. To maintain a semblance of brevity, these cases are excluded here except by reference to the preliminary version.

---

[*] mike@geometrics-unlimited.com

It is important to recognize that the identification of particular results from comparison with a tabulated database is a powerful adjunct to the WZ method, since that method can usually be counted on to generate a recursion formula when an integer parameter *"n"* exists. Without some means of obtaining the starting values for the recursion, or solving the recursion itself[3], it is of little practical utility. As will be demonstrated later, the conjunction of the two methods gives a powerful means of evaluating generalized representations of Dixon's, Watson's and Whipple's theorem using a computer algebra program.

**2. Method**

Prudnikov et. al.[7], Section 7.4.4 tabulate many closed form relations for special parametric cases of $_3F_2(\ldots|1)$. Each of the 10 Thomae relations listed in Appendix A was applied to 70[*] of the closed-form identities listed there, resulting in 630 possibly new and/or different closed-form sums. Each of these results was compared against all the others, by searching for the existence of a valid transformation among all the top and bottom parameters taking account of the symmetry that exists among these parameters. If such a transformation was found, the relation that had the most general parametric set was retained, and the particular (equivalent or special) case was discarded. Additionally, all results having a parametric excess equal to a non-positive integer were discarded, unless that particular case corresponded to a terminating sum, as were all cases (Karlsson-Minton[9,10,11]) where a top parameter exceeded a bottom parameter by a small positive integer. Finally, all terminating sums that reduce to other terminating sums with the same number of terms were removed. This procedure yielded a "base set" of fundamental results obtainable from Ref. 7. The list of input cases was then expanded to include other results harvested from the literature for which an acceptable transformation from the existing "base set" parameters could not be found, and the same procedure was applied.

There were 99 survivors of this culling procedure, of which 23 covered the set of evaluations initially input, leaving 76 relations that are thought to be new (but see Section 3). Of these, 42 involve transformations that reduce an infinite series to a finite series – these are listed elsewhere[8]. The remaining 35 "*closed-form*"[†] results are listed in Appendix B, and, with one exception, all were included in I. Recent new work by Maier has expanded the basic set of new results from which candidates could be harvested, and this has yielded another 6 sums, as well as another 27 sets of relations related by recursion. All of the entries so-found were tested for novelty using computer algebra simplification commands and for validity using numerical evaluation. In none of these cases was a computer algebra code able to arrive at the closed-form result given in the appendices.

---

[*] It is recognized that many of these results were special cases of more general parameter sets.
[†] *"closed form"*: does not involve a specific (infinite) summation operation

**3. Input cases**

*3.1 Notation*

Throughout, I use the symbols *"n"* and *"m"* to represent positive integers giving a strong constraint on the acceptable transformations among parameters. Otherwise all symbols represent arbitrary complex, continuous variables. In many of the formulae quoted, it is recognized that the real part of the *"parametric excess"* ($b_1 + b_2 - a_1 - a_2 - a_3$) of a particular hypergeometric function on the left-hand side must exceed zero in order for an infinite series representation to exist. However, the right-hand side of any such equation is a valid representation of that hypergeometric function by the principle of analytic continuation. Therefore, provided that the parametric excess is not a negative, absolute constant, any such limitation is irrelevant and is only specifically quoted as a convenient categorization parameter (*"excess"*).

*3.2 Comments on Section 7.4.4 of Ref. 7.*

Each of the formula copied from Section 7.4.4 was tested numerically to guard against transcription errors. In this way, it was discovered that a number of formulae given in that source are incorrect. In particular, [7.4.4.19] does not satisfy any numerical test for arbitrary values of the parameter *"b"*, except in the case that $b = -n$. Thus this equation was retained (with an obvious change of notation) in the form

$$_3F_2\left(\begin{matrix} a, b, -n \\ b+\tfrac{1}{2}, a-n+\tfrac{1}{2} \end{matrix} \middle| 1 \right) = \sqrt{\pi}\, \frac{\Gamma(a-n+\tfrac{1}{2})\Gamma(b+\tfrac{1}{2})\Gamma(b-a+n+\tfrac{1}{2})}{\Gamma(a+\tfrac{1}{2})\Gamma(-n+\tfrac{1}{2})\Gamma(b-a+\tfrac{1}{2})\Gamma(b+n+\tfrac{1}{2})} \qquad (1)$$

Equations [7.4.4.38] and [7.4.4.73] do not satisfy any non-trivial numerical tests, and lacking a source reference from which a corrected form might have been obtained, these two results were omitted from consideration.

Several misprints were also discovered. Notably, equation [7.4.4.43] should read

$$_3F_2\left(\begin{matrix} 1, 2, a \\ 3, b \end{matrix} \middle| 1 \right) = -2(b-2) + 2(b-1)^2 \psi'(b) \qquad (2)$$

equation [7.4.4.55] should read

$$_3F_2\left(\begin{array}{c}a, a-n, a-n\\a-n+1, a-n+1\end{array}\bigg|1\right) = \frac{\pi(n-1)!(n-a)}{(1-a)_{n-1}\sin a\pi} \tag{3}$$

equation [7.4.4.67] becomes

$$_3F_2\left(\begin{array}{c}2, a, a\\a+2, a+2\end{array}\bigg|1\right) = -a^2(a+1)^2[(2a-1)\psi'(a)-2] \tag{4}$$

and equation [7.4.4.71] should be

$$_3F_2\left(\begin{array}{c}3, a, a+1\\a+2, a+3\end{array}\bigg|1\right) = (1+a)^2(2+a)(-2a^2(a-1)\psi(1,a)+a(2a-1))/4. \tag{5}$$

Based on [Ref. 7, 7.4.4.25] relating a particular $_3F_2(1)$ to a special $_2F_1(-1)$ three further results were developed. These are referred to as *"Prudnikov 7.4.4.25 variation 1-3"* in Ref. 8 and are respectively given below:

$$_3F_2\left(\begin{array}{c}a, a+\frac{1}{2}, b\\b+\frac{3}{2}-a, b-a+1\end{array}\bigg|1\right) = \frac{2^{-b}\pi^{\frac{1}{2}}\Gamma(2+2b-2a)\Gamma(2+b-4a)}{\Gamma(2+2b-4a)(2a-1)}$$
$$\times\left(\frac{1}{\Gamma(\frac{1}{2}b)\Gamma(\frac{3}{2}+\frac{1}{2}b-2a)} - \frac{1}{\Gamma(\frac{1}{2}b+\frac{1}{2})\Gamma(1+\frac{1}{2}b-2a)}\right) \tag{6}$$

$$_3F_2\left(\begin{array}{c}a, a+\frac{1}{2}, 1\\\frac{3}{2}-a-\frac{1}{2}n, 2-a-\frac{1}{2}n\end{array}\bigg|1\right) = \frac{\Gamma(3-n-2a)\Gamma(2-n-4a)}{\Gamma(3-n-4a)}$$
$$\times\left(\frac{2^{-4a-n}\Gamma(1-2a)}{\Gamma(2-n-4a)} + \frac{1}{2\,\Gamma(2a)}\sum_{k=0}^{n-1}\frac{(-1)^k\Gamma(2a+k)}{\Gamma(2+k-n-2a)}\right) \tag{7}$$

$$_3F_2\left(\begin{array}{c}a, a+\frac{1}{2}, 1\\\frac{3}{2}-a+\frac{1}{2}n, 1-a+\frac{1}{2}n\end{array}\bigg|1\right) = \Gamma(1-2a)\frac{(1+n-2a)}{(1+n-4a)}$$
$$\times\left(\frac{2^{n-1-4a}\Gamma(1+n-2a)}{\Gamma(1+n-4a)} - \frac{1}{2\,\Gamma(2a-n)}\sum_{k=1}^{n-1}\frac{(-1)^k\Gamma(2a-n+k)}{\Gamma(1+k-2a)}\right) \tag{8}$$

*3.3 Other results extracted from the literature*

Two general results given elsewhere[12] survived the tests described in the previous section. They are quoted below and referred to respectively as *"Ref. 12, Lemma 2.1"* and *"Ref. 12, Lemma 2.2"*:

$$_3F_2\left(\begin{matrix}1,n+1,a\\n+2,b\end{matrix}\Big|1\right) = \frac{(n+1)\Gamma(b)}{\Gamma(a)}\left(\frac{\Gamma(a-n-1)(\psi(b-n-1)-\psi(b-a))}{\Gamma(b-n-1)}\right.$$
$$\left. - \sum_{l=0}^{n-1}\frac{\Gamma(a+l-n)}{\Gamma(b+l-n)(l+1)}\right) \quad (9)$$

and

$$_3F_2\left(\begin{matrix}a,b,c\\n+b,c+1,\end{matrix}\Big|1\right) = \frac{(b)_n\Gamma(c+1)\Gamma(1-a)}{(b-c)_n\Gamma(c+1-a)} +$$
$$c\Gamma(b+n)\Gamma(c-b+1-n)\sum_{l=0}^{n-1}\frac{\Gamma(n-l-a)(-1)^l}{\Gamma(b+n-a-l)\Gamma(n-l)\Gamma(c-b-n+2+l)} \quad (10)$$

Notice that a result of Rao et al.[13] citing Ramanujan's notebook corresponds to the special case $b=c, n\to n+1$ of **(10)**[14].

Several years ago, Sharma[15] claimed to have found two new closed forms for a particular $_3F_2(\ldots|1)$, based on special choices of parameters for a new $_4F_3(\ldots|1)$ obtained by evaluating a double series. Unfortunately, the new $_4F_3(\ldots|1)$ given in Sharma's paper does not satisfy any numerical tests, and is clearly incorrect, since the right hand-side is symmetric in (Sharma's) variables $\alpha$ and $\rho$, while the left-hand side is not. Using the transformation algorithm described previously, it is noted that the left-hand side of Sharma's equation (7) is a special case of [Ref. 7, 7.4.4.25] and the left-hand side of Sharma's equation (8) corresponds to a special case of [Ref. 7, 7.4.4.47], which in turn is a special case of Watson's theorem $W_{0,0}$ (see Section 3.4). However, the right-hand sides of Sharma's equations do not correspond to the respective right-hand sides of the quoted results of Ref. 7. For these reasons, Sharma's results were omitted from consideration.

Recently, Exton[16] gave a "new" result for a special case of $_3F_2(\ldots|1)$ based on a purportedly "new" two-term transformation[17] between different $_3F_2(\ldots|1)$. It is easily seen that Exton's "new" transformation, is actually a symmetric permutation of **T3** of Appendix A, and that the "new" result is incorrect[18] by calculating the different limits *a, b* or *c=0* on both sides of Exton's equation (13). For this reason, that result was also omitted. Exton has also given some new identities[19] and sums for $_3F_2(\ldots|1)$ by evaluating one part of a double series. Each of these results [Ref. 19, 1.13 and 1.14] corresponds to a special case of results given in Appendix B. Additionally, the case [3.3, q=2] of another of Exton's new results[20] is a special case of **(Ref. 8, Eq. B.1)**.

All of the relevant results given in Table 2 of Krupnikov and Kölbig[21] are special cases of results found here; however two new results found in a paper of Gessel and Stanton[22] – equations (1.6) and (1.9) respectively - are listed below:

$$_3F_2\left(\begin{matrix} 2a, 1-a, -n \\ 2a+2, -a-\frac{1}{2}-\frac{3}{2}n \end{matrix} \Big| 1\right) = \frac{((n+3)/2)_n (n+1)(2a+1)}{(1+\frac{1}{2}(n+2a+1))_n (2a+n+1)} \tag{11}$$

$$_3F_2\left(\begin{matrix} -sa+s+1, a-1, -n \\ a+1, -s(n+a)-n \end{matrix} \Big| 1\right) = \frac{a(1+s+sn)_n (n+1)}{(1+s(a+n))_n (a+n)} \tag{12}$$

The pedigree of **(12)** is interesting – Gessel and Stanton[*] only give a proof for $n>0$, citing an unpublished letter from Gosper to Askey for the case $n<0$. Since **(12)** satisfies numerical tests for $n<-1$, a second version of that result with $n \to -n-1$ was included in the database. Note that **(12)** cannot be certified by the WZ method[3]. See also Section 4.

In addition, Ref. 22 also gives, in equation (5.16), a special case of a $_4F_3$, citing a problem posed by Fields and Luke. Since this equation is of Minton type (an upper parameter exceeds a lower parameter by unity), it can be reduced as follows:

$$_4F_3\left(\begin{matrix} 1, 1+a/(s-1), bs-a, -n \\ a/(s-1), b+1, 1-a-ns \end{matrix} \Big| 1\right) = {}_3F_2\left(\begin{matrix} 1, bs-a, -n \\ b+1, 1-a-ns \end{matrix} \Big| 1\right)$$
$$- \frac{(s-1)(sb-a)n}{a(b+1)(1-a-ns)} {}_3F_2\left(\begin{matrix} 2, 1+bs-a, 1-n \\ b+2, 2-a-ns \end{matrix} \Big| 1\right) \tag{13}$$

The first of the two terms on the right can be summed (Ref. 7, 7.4.4.13) leading to:

$$_3F_2\left(\begin{matrix} 2, 1+sb-a, 1-n \\ b+2, 2-a-ns \end{matrix} \Big| 1\right) = \frac{(b+1)(a-1+ns)(a+ns)ab\Gamma(n)\Gamma(1-ns+n+b-sb)}{(sb-a)(s-1)\Gamma(b+n+1)\Gamma(1-a-ns+n)}$$
$$\times \sum_{L=0}^{n} \frac{\Gamma(-a-ns+L)\Gamma(b+L)}{\Gamma(L+1)\Gamma(1-ns+b-sb+L)} \tag{14}$$
$$+ \frac{b(a+ns)(-1+a+ns)(b+1)}{n(s-1)(sb-a)(b+n)}$$

Although this result does not appear to be very profound – it simply transforms a finite sum – when **T1-T9** are applied to this equation a number of new results emerge that reduce infinite series of Minton type into finite sums. These are listed in Appendix B of Ref. 8; the genesis of each result is labeled *"Gessel and Stanton, (Eq. 5.16)"*. Notice that (Ref. 8, **B.38**) based on **(14)** generalizes (Ref. 8, **B.53**), which is based on **(12)**. As well, **T6** applied to **(14)** reduces to **(11)** as a special case.

---

[*] In Ref. 22, Eq. 5.16, which sums a terminating $_4F_3$, tested numerically true for continuous "*n*". As well, "strange" Eqs. 3.11 through 3.16 all reduce to a $_2F_1$ so they were omitted.

**3.4** *Generalized (Contiguous) Watson, Dixon and Whipple's theorem*

Lavoie[23] has given two "probably new" summation formulae contiguous to Watson's theorem. These were checked with the database; it was discovered that Lavoie's equation (2) is new[*], and it was therefore added; Lavoie's "new" equation (1) can be obtained from equation (2) of that same paper by the application of Thomae relation **T1** so it was excluded. All other relations quoted in that paper are special cases of results already included in the database[†] as are all the $_3F_2(\ldots|1)$ relations recorded in another of Lavoie's works[24]. For completeness' sake, Lavoie's equation (2) is reproduced below:

$$_3F_2\left(\begin{array}{c}a,b,c\\ \tfrac{1}{2}(a+b+1), 2c+1\end{array}\,\middle|\,1\right) = \frac{2^{a+b}\Gamma(\tfrac{1}{2}a+\tfrac{1}{2}b+\tfrac{1}{2})\Gamma(c+\tfrac{1}{2})\Gamma(c-\tfrac{1}{2}a-\tfrac{1}{2}b+\tfrac{1}{2})}{\Gamma(\tfrac{1}{2})\Gamma(a+1)\Gamma(b+1)}$$
$$\times\left(\frac{\Gamma(\tfrac{1}{2}a+1)\Gamma(\tfrac{1}{2}b+1)}{\Gamma(c-\tfrac{1}{2}a+\tfrac{1}{2})\Gamma(c-\tfrac{1}{2}b+\tfrac{1}{2})} - \frac{ab\,\Gamma(\tfrac{1}{2}a+\tfrac{1}{2})\Gamma(\tfrac{1}{2}b+\tfrac{1}{2})}{4\Gamma(c-\tfrac{1}{2}a+1)\Gamma(c-\tfrac{1}{2}b+1)}\right) \quad (15)$$

Additionally, in several other works Lavoie et. al. give new results for cases contiguous to Dixon's,[25] Watson's[26] and Whipple's[27] theorems. As shown by Bailey[28] in the case of a $_3F_2$, Watson's theorem can be easily obtained by the application of (a symmetric permutation of) **T2** of Appendix A to Dixon's theorem, and Whipple's theorem can be subsequently obtained by the application of **T2** to Watson's theorem. It turns out that these three theorems are similarly closed under the application of any of the other Thomae relations, as are the contiguity identities for Whipple's, Watson's and Dixon's theorems studied by Lavoie et. al. Since the results given by Lavoie et. al. for Whipple's and Dixon's theorems are limited to a small subset of near-diagonal contiguous cases[‡], whereas Lewanowicz[29] has independently given a general result valid for all off-diagonal instances of the generalized Watson's theorem, in this sense Watson's theorem may be more fundamental than the other two. Therefore, in principle, only Lewanowicz' general result for contiguous elements of the generalized Watson's theorem needs to be retained in any database, with the other relations obtained by the application of any of the Thomae transformations (see below).

However, the complexity of these results and the embodiment of the WZ method in computer algebra codes leads to another practical method of obtaining any particular element of the generalized result for any of these cases. Consider the elements $X_{m,n}$ contiguous to the (well-poised) Dixon's element ($m=n=0$):

$$X(a,b,c,m,n) = X_{m,n} \equiv {}_3F_2\left({a,b,c \atop 1+m+a-b, 1+m+n+a-c} \mid 1\right) \tag{16}$$

with

$$X(a,b,c,n,m) = X(a,c,b,m+n,-n) \tag{16b}$$

---

[*] A correspondent has pointed out that "*this result and a similar result were already obtained by B.N. Bose in his paper entitled "On some transformations of the generalized hypergeometric Series", published in the J. Indian Math. Soc., Vol. VIII, nos. 3 and 4, 119-128 (year unspecified). This result is given in that paper as equation (17), page 128*". I have been unable to obtain that paper.
[†] Is it possible that Lavoie's results are the source for several entries in Prudnikov's table?
[‡] 38 cases for Whipple's theorem, 39 for Dixon's theorem.

Computer algebra codes will yield (very lengthy) three part recursion formulae in the variables *m* and *n*. With knowledge of four values of $X_{m,n}$ and using the recursion so-obtained*, it is then possible to obtain any of the sums contiguous to Dixon's theorem ($X_{0,0}$) for any value of *m* and *n*. The case $X_{0,0}$ is well known and given by [Ref. 7, 7.4.4.21], $X_{0,1}$ is given by **(B.14)**, $X_{0,-1}$ corresponds to [Ref. 7, 7.7.4.20], $X_{1,0}$ is given by [Ref. 7, 7.7.4.22] and $X_{1,-1}$ is given (symmetrically) by **(B.14)**.

A similar method can be used to obtain generalizations of Watson's theorem, where the recursion is far simpler. It is found that any of the elements $W_{m,n}$ defined by the generalized Watson theorem

$$W_{m,n} \equiv {}_3F_2\left({a,b,c \atop \tfrac{1}{2}(a+b+1+m), 2c+n} \mid 1\right) \tag{17}$$

obey three-part recursion formulae in "*m*" and "*n*" given by **(C.1)** and **(C.2)** so it remains necessary to locate eight closed forms to start the recursion for all values of "*m*" and "*n*". The element $W_{0,0}$ corresponding to Watson's theorem is well-known and given by [Ref. 7, 7.4.4.18] or [Ref. 30, 3.13.3(7)]; $W_{0,1}, W_{0,-1}, W_{-1,1}, W_{1,0}$ and $W_{1,-1}$ are given by **(15)**, **(B.7)**, **(B.8)**, **(B.11)** and **(B.15)** respectively. Other starting results for the initial elements $X_{m,n}$ can be obtained from Lavoie et. al.[25], and corresponding elements $W_{m,n}$ follow from **(21)** (see below). As a test, the elements $X_{-1,0}$ and $X_{2,0}$, obtained by recursion (listed in **(C.3)** and **(C.4)**), were compared to, and agreed with, the same elements taken from Ref. 25. Using **(21)**, the corresponding elements $W_{-1,0}$ and $W_{2,0}$ were found; these are given by **(C.5)** and **(C.6)** respectively. With this method, any of the other elements $X_{m,n}$ or $W_{m,n}$ can now be found by recursion in either $X_{m,n}$ or (preferably) $W_{m,n}$ using **(C.1)** and **(C.2)**.

---

* Note that Lewanowicz' coefficients are also given recursively.

The third member of the trio consists of generalized elements contiguous to Whipples' theorem. As noted by Lavoie et. al.[27] using a method similar to that employed here, elements $P_{m,n}$ of Whipple's theorem

$$P_{m,n} \equiv {}_3F_2\left(\begin{array}{c} a,\ b, 1-b+m+n \\ c,\ 1+2a+m-c \end{array} \Big| 1\right) \tag{18}$$

are related to $X_{m,n}$. The relationship can be found by applying **T8** to $X_{m,n}$ giving

$$P_{m,n} = \frac{X(2a-b-c+m+1, 1+a-c+m, 1-b+m+n, m, n)\Gamma(a-n)\Gamma(c)}{\Gamma(b+c-1-m-n)\Gamma(a-b+m+1)} \tag{19}$$

The correspondence with Watson's theorem is obtained by applying **T1** to $W(a,b,c,m,n) \equiv W_{m,n}$ giving

$$P_{m,n} = \frac{W(c-b, 2a-c+m+1-b, a-n, m, n)\Gamma(a-n)\Gamma(c)\Gamma(1+m+2a-c)}{\Gamma(b)\Gamma(2a-n)\Gamma(a-b+m+1)} \tag{20}$$

Notice that **(19)** fails[30] for the case $m=n=0$ when "$a$" is a non-positive integer, unless "$b$" is an integer[*].

As an adjunct to the above, apply **T2** to $X_{m,n}$ to find

$$\begin{aligned}W_{m,n} = X(2c+n-a, (b-a+m+1)/2, (1+m-a-b)/2+c+n, m, n) \\ \times \frac{\Gamma(c+n+\tfrac{1}{2}(1+m-a-b))\Gamma(2c+n)\Gamma(\tfrac{1}{2}(a+b+1+m))}{\Gamma(a)\Gamma(c+n+\tfrac{1}{2}(1+m+b-a))\Gamma(2c+n+\tfrac{1}{2}(1+m-a-b))}\end{aligned} \tag{21}$$

or

$$\begin{aligned}X_{m,n} = W(1+m+a-2b,\ a,\ 1+m+a-b-c,\ m,\ n) \\ \times \frac{\Gamma(a-2b-2c+2+2m+n)\Gamma(1+a-c+m+n)}{\Gamma(1-c+m+n)\Gamma(2a-2b-2c+2+2m+n)}\end{aligned} \tag{22}$$

See also **(B.12)** and **(B.13)**.

*3.5 Other computer applications*

Although others have also made extensive use of computer algebra to obtain transformations and identities between generalized hypergeometric series, I found no other new sums in the literature when the preliminary version of this work was prepared (see Section 5 below). In particular, Gessel[31] gives a large number of new identities of which only one (equation 30.1) is relevant to this work. That result is of the Karlsson-Minton type being a special case of [Ref. 7, 7.4.4.10].

---

[*] see ref. 30 Eqs 3.13(8) and (9). This exceptional case is given by Prudnikov 7.4.4.104.

Krattenthaler's Mathematica package HYP[32] and Gauthier's Maple package HYPERG,[33] contain a number of equivalent rules for evaluating $_3F_2(…|1)$ - HYP taken from the limiting case ($q \to 1$) of basic hypergeometric series[34], HYPERG taken from Slater[35]. Although not obvious, several of these are special cases of Appendix B. The correspondence follows (see Table 1, below):

**Table 1. Showing the correspondence between the database sums of refs. 32, 33, 35 and the database elements described herein.**

| Equation reference | database number and equivalent special case |
|---|---|
| S3201 | **Ref. 8, B.3** with $m=1$ |
| S3204 | [Ref. 7, 7.4.4.10] with $e = 1 - b + \lambda$; $c = 1 + \lambda/2$ |
| S3235 | **Ref. 8, B.3** with $m = 2$; $c = 1 + a - b$ |
| Slater III.15 | [Ref. 7, 7.4.4.10] with $c = 1 + \lambda/2$; $b = -n$ |

### 4. Comments and a Conjecture Resolved

As discussed, Appendix B gives new closed-form results for a number of hypergeometric sums. Other results, based on the two-part Thomae relations, contain a series that is also equivalent to the well-known three part relations between hypergeometric $_3F_2(1)$[30], but, usefully, still reduce an infinite series with a parameter *"n"* to a finite series with *"n"* (or less) terms – see Ref. 8 for a listing of such cases. In a number of instances, new results were obtained which were deleted from Appendix B because they trivially reduced to Karlsson-Minton type.

In a recent work, Krattenthaler and Rivoal[36] introduced two new three-parameter two-part relations between $_3F_2(1)$. These were tested against the database, and 31 results were found that satisfied one or the other of these relations. However, when the appropriate equality was applied to any of these cases in the same manner as was done for the Thomae relations, the resulting sum was found to be already included in the database, so no new sums were found from this procedure.

In the course of testing the results, an observation was made that leads to two new sums and corresponding Thomae progeny. It was noticed that **(11)**, due to Gessel and Stanton[22], satisfies numerical tests for $n \to -2n$, leading to the following new result:

$$_3F_2\left(\begin{matrix}2a, 1-a, 2n \\ 2a+2, -a+3n-\frac{1}{2}\end{matrix}\Big|1\right) = \frac{(2a+1)\Gamma(a-n+\frac{1}{2})\Gamma(-3n+\frac{5}{2})}{3\Gamma(\frac{3}{2}-3n+a)\Gamma(\frac{3}{2}-n)} \quad n \geq 0 \quad \textbf{(23)}$$

which was verified using the WZ certification procedure[3]. If $n = 1$, **(23)** reduces to **(B.20)** (with the substitutions $a \to 2a$, $b \to 1-a$ ).

In [8], it was observed that [Ref. 22, (5.16)] survived all numerical tests when *"n"* was not integral. Setting $a = s - 1$ in that equation reduces it to a $_3F_2$, which turns out to be equivalent to **(B.18)**. Replace *"n"* by a continuous variable and with a change of notation in [8], I conjectured that:

$$_3F_2\left(\begin{matrix} 2, a, b \\ c, (2c - 3 - a + ab - b)/(c - 2) \end{matrix} \middle| 1\right) = -\frac{(c-1)(c + ab - b - a - 1)}{(c - 1 - b)(a + 1 - c)} \quad (24)$$

The application of T1 to **(24)** yielded another (conjectured) new result – see **(B.34),** and, due to an error in my processing code, the application of T8 missed another, which is now given here by **(B.35)**. It has recently come to my attention that all three of these results are summarized as **Theorem 3.4** in a paper by Maier [37] who indicates that **(B.35)** was first given by Slater ([35], (2.1.1.10)), but misprinted. As well, Chu [38] has independently derived **(24)** by a method of creative telescoping. Thus **(24)**, **(B.34)** and **(B.35)** now have solid pedigrees.

**5. Extensions**

Working from [37], Maier has broadened his results to more general hypergeometric series [39], which, in passing, yield several new cases of $_3F_2(\ldots|1)$. From ([39], (1.2)) we find

$$_3F_2\left(\begin{matrix} a, a + \tfrac{1}{3}, a + \tfrac{2}{3} \\ 3a + 1, \tfrac{1}{2} \end{matrix} \middle| 1\right) = \frac{3^{3a}(1 + 4^{-3a})}{2} \quad (25)$$

from which nine new Thomae progeny are obtained with no counterpart in the previous database. From that same paper, Theorem (7.1) with Maier's index $\ell = 0$ and $\ell = -1$ respectively, yields

$$_3F_2\left(\begin{matrix} 2a, 2a - \tfrac{1}{3}, 2a + \tfrac{1}{3} \\ 3a, 3a + \tfrac{1}{2} \end{matrix} \middle| 1\right) = \left(\tfrac{3}{2}\right)^{6a-1} \quad (26)$$

and

$$_3F_2\left(\begin{matrix} 2a, 2a - \tfrac{1}{3}, a + \tfrac{1}{3} \\ 1 + 3a, 3a + \tfrac{1}{2} \end{matrix} \middle| 1\right) = \left(\tfrac{3}{2}\right)^{6a-1} / (2a + \tfrac{2}{3}) \quad (27)$$

Again, nine new Thomae progeny are found from each of these. In general, any further values corresponding to Maier's index $\ell = n$ can be obtained by recursion. For these more general values of $\ell$, the simplest (of some complicated) three-part recursion formulae in both the forward[*] and backward directions, belong to the second Thomae progeny at $\ell = n$, that is

$$_3F_2\left(\begin{array}{c}½-n, a-n/2, a+½-n/2\\ 2a+1/6-n, 2a+5/6-n\end{array}\bigg|1\right) = V_2(a,n) \tag{28}$$

where

$$V_2(a,0) = \frac{3^{(6a-1)}\Gamma(2a)\Gamma(1/6+2a)\Gamma(2a+5/6)}{\sqrt{\pi}\Gamma(6a)}$$

$$V_2(a,-1) = 4\frac{3^{(6a-1)}\Gamma(2a)\Gamma(7/6+2a)\Gamma(2a+11/6)}{\pi\Gamma(1+6a)(2a+2/3)} \tag{29}$$

The corresponding backward/forward recursion formula for $V_2(a,n)$ are given in Appendix C by **(C.7)** and **(C.8)**. The nine Thomae progeny of **(28)** are listed in Appendix D, equations numbers **(D.1)** to **(D.9)**. In that listing, **(D.5), (**being the fifth Thomae progeny of its own second progeny), corresponds to the generalization of the original cases **(26)** and **(27)** to $\ell = n, n > 0$. The labeling follows the notation of Appendix B, with the addition of the "**:T2**" symbol to remind that each case was derived by applying the respective transformation to the **T2** progeny of **(26)** and **(27)**.

From Theorem (7.3) of [39] corresponding to Maier's indices $\ell = 0, \kappa = 1$ and **(25)** which corresponds to the same theorem with $\ell = 0, \kappa = 0$, another set of new cases of $_3F_2(\ldots|1)$ exists, being:

$$_3F_2\left(\begin{array}{c}a, a+1/3, a+2/3\\ 3a, 3/2\end{array}\bigg|1\right) = -(3)^{3a}(8^{1-2a}-1)/(12a-6) \tag{30}$$

$$_3F_2\left(\begin{array}{c}a, a+1/3, a+2/3\\ 3a+2, 1/2\end{array}\bigg|1\right) = (3)^{3a}(8^{-2a}(9a+4)+4)/(16a+8) \tag{31}$$

There are two possible generalizations, the first based on recursions between **(30)** and **(25)**, the second between **(31)** and **(25)**. Define

---

[*] For the original case, corresponding to **(26)** and **(27)**, recursion fails in the forward (i.e. $n > 0$) direction, since the *parametric excess* is ½-n<0.

$$_3F_2\left({a, a+\tfrac{1}{3}, a+\tfrac{2}{3} \atop 3a-n, \tfrac{1}{2}+n} \mid 1\right) \equiv U(a,n) \tag{32}$$

and

$$_3F_2\left({a, a+\tfrac{1}{3}, a+\tfrac{2}{3} \atop 3a+n, \tfrac{1}{2}} \mid 1\right) \equiv V(a,n) \tag{33}$$

The nine new Thomae progeny derived from **(32)** are labeled "Maier Eq. (7.3) with L=n, Variation 1" and numbered **(D.10)** to **(D.18)** in Appendix D; those derived from **(33)** are labeled similarly as "…Variation 2" and numbered **(D.19)** to **(D.27)**. **(C.9)** and **(C.10)** list the forward /backward recursion formula for $U(a,n)$; the right-hand sides of **(25)** and **(30)** defining the starting values $U(a,0)$ and $U(a,1)$ respectively. **(C.10)** lists the forward recursion formula for $V(a,n)$; the right-hand sides of **(25)** and **(31)** define the starting values $V(a,1)$ and $V(a,2)$ respectively (**(33)** diverges in the backward direction $n \leq 0$ ).

## 6. Summary

Sixty-two new sums for hypergeometric $_3F_2(…|1)$ were found and a large number of relations that reduce an infinite series to a finite series were reported elsewhere. This review fairly well summarizes the existing literature to the best of my knowledge. Together with previously known forms, this collection forms a useful basic database for the possible computerized identification of any desired sum of the form $_3F_2(…|1)$, by simply seeking a transformation between the parameters of the candidate and elements of the database.  The existence of such a database is an important adjunct to the WZ recursion method of evaluating such sums, since it provides a starting point for the recursion in several important cases.  The validity of two forms previously conjectured was verified and a number of errors appearing in the literature were noted.

Appendix A

The following lists 10 inequivalent Thomae relations for $_3F_2([a, b, c], [f, e], 1)$ as they are numbered and referenced in the text. There are 110 other symmetric permutations. T10 is the identity.

T1:

$$_3F_2([e+f-a-b-c, f-c, e-c], [e-b+f-c, e+f-a-c], 1)$$
$$\times \Gamma(e+f-a-b-c)\Gamma(f)\Gamma(e)/(\Gamma(c)\Gamma(e-b+f-c)\Gamma(e+f-a-c))$$

T2:

$$_3F_2([e+f-a-b-c, -b+f, e-b], [e-b+f-c, e-b+f-a], 1)$$
$$\times \Gamma(e+f-a-b-c)\Gamma(f)\Gamma(e)/(\Gamma(b)\Gamma(e-b+f-c)\Gamma(e-b+f-a))$$

T3:

$$\frac{_3F_2([f-c, -b+f, a], [e-b+f-c, f], 1)\,\Gamma(e+f-a-b-c)\,\Gamma(e)}{\Gamma(e-a)\,\Gamma(e-b+f-c)}$$

T4:

$$\frac{_3F_2([e-c, e-b, a], [e-b+f-c, e], 1)\,\Gamma(e+f-a-b-c)\,\Gamma(f)}{\Gamma(f-a)\,\Gamma(e-b+f-c)}$$

T5:

$$_3F_2([e+f-a-b-c, f-a, e-a], [e+f-a-c, e-b+f-a], 1)$$
$$\times \Gamma(e+f-a-b-c)\Gamma(f)\Gamma(e)/(\Gamma(a)\Gamma(e+f-a-c)\Gamma(e-b+f-a))$$

T6:

$$\frac{_3F_2([f-c, f-a, b], [e+f-a-c, f], 1)\,\Gamma(e+f-a-b-c)\,\Gamma(e)}{\Gamma(e-b)\,\Gamma(e+f-a-c)}$$

T7:

$$\frac{_3F_2([e-c, e-a, b], [e+f-a-c, e], 1)\,\Gamma(e+f-a-b-c)\,\Gamma(f)}{\Gamma(-b+f)\,\Gamma(e+f-a-c)}$$

T8:

$$\frac{_3F_2([-b+f, f-a, c], [e-b+f-a, f], 1)\,\Gamma(e+f-a-b-c)\,\Gamma(e)}{\Gamma(e-c)\,\Gamma(e-b+f-a)}$$

T9:

$$\frac{{}_3\mathrm{F}_2([e-b,\ e-a,\ c],\ [e-b+f-a,\ e],\ 1)\,\Gamma(e+f-a-b-c)\,\Gamma(f)}{\Gamma(f-c)\,\Gamma(e-b+f-a)}$$

$T10:\qquad {}_3\mathrm{F}_2([a,\ b,\ c],\ [f,\ e],\ 1)$

## Appendix B

The following lists 35 new evaluations of ${}_3\mathrm{F}_2([a,\ b,\ c],\ [e,\ f],\ 1)$ labelled B.1 to B.35, obtained as discussed in the text. Each is classified by the (parametric) "excess". The derivation can be duplicated by applying the appropriate Thomae relation (Appendix A) to the basic relation whose source is given. For example, to obtain B.3, apply T3 of Appendix A to Eq. 7.4.4.17 of Ref. 7.

$$\text{excess} = a$$

**B.1  Prudnikov 7.4.4.22  : T7**

$${}_3\mathrm{F}_2([a,\ b,\ -b+2],\ [c,\ 2a+2-c],\ 1)\ =$$

$$\left(\left(\frac{1}{2}\frac{\Gamma(a-\tfrac{1}{2}c+1-\tfrac{1}{2}b)\,\Gamma(\tfrac{1}{2}b-1+\tfrac{1}{2}c)}{\Gamma(a+\tfrac{1}{2}b-\tfrac{1}{2}c)\,\Gamma(-\tfrac{1}{2}b+\tfrac{1}{2}c)} - \frac{1}{2}\frac{\Gamma(\tfrac{3}{2}+a-\tfrac{1}{2}b-\tfrac{1}{2}c)\,\Gamma(\tfrac{1}{2}b+\tfrac{1}{2}c-\tfrac{1}{2})}{\Gamma(a-\tfrac{1}{2}c+\tfrac{1}{2}+\tfrac{1}{2}b)\,\Gamma(-\tfrac{1}{2}b+\tfrac{1}{2}c+\tfrac{1}{2})}\right)\right.$$

$$\left.\Gamma(2a+2-c)\,\Gamma(c)\right)\Big/\left((-b+1)\,(-c+a+1)\,\Gamma(c+b-2)\,\Gamma(2+2a-b-c)\right)$$

$$\text{excess} = -b+c-a+2$$

**B.2  Prudnikov 7.4.4.17  : T2**

$${}_3\mathrm{F}_2([a,\ 2,\ b],\ [c,\ 4],\ 1)\ =$$

$$-6\frac{(2c-5+b-ab+a)\,\Gamma(c)\,\Gamma(-b+2+c-a)}{(a-3)(b-1)(-2+a)(a-1)(b-3)(b-2)\,\Gamma(c-a)\,\Gamma(c-b)}$$
$$+\frac{6(c-2)(c-1)(ab-3b-3a+3+2c)}{(a-3)(b-1)(-2+a)(a-1)(b-3)(b-2)}$$

$$\text{excess} = 2$$

**B.3  Prudnikov 7.4.4.17  : T3**

$${}_3\mathrm{F}_2([a,\ 2,\ b],\ [c,\ 4+a-c+b],\ 1)\ =$$

$$\frac{(-a-2+c-b)(-a-3+c-b)(c-2)(c-1)(3-4c-ca+c^2+3a+ab+3b-bc)}{(-a-2+c)(c-3-b)(c-2-b)(c-1-b)(-a-3+c)(-a+c-1)}$$
$$+\frac{(a+5-ca+c^2-4c+ab-bc+b)\,\Gamma(4+a-c+b)\,\Gamma(c)}{(-a-2+c)(c-3-b)(c-2-b)(c-1-b)(-a-3+c)(-a+c-1)\,\Gamma(b)\,\Gamma(a)}$$

excess= 2

B.4  Gessel & Stanton, SIAM J. Math. Anal. 13,2,295(1982)  Eq(1.9)
with n->-1-n    : T7

$${}_3F_2([a,\, b,\, -\frac{n\,(n-b+a)}{b-n}],\, [1+b-\frac{a\,n}{b-n},\, 1+a+n],\, 1) \;=\;$$

$$\frac{(-1)^n \sin(\pi b)\,\Gamma(1+b+\dfrac{a\,n}{-b+n})\,\Gamma(-b)}{\sin(\pi a)\,\Gamma(1+a)\,\Gamma(-a-n)\,\Gamma(n+1+\dfrac{a\,n}{-b+n})}$$

excess = 2

B.5   Gessel and Stanton with n->-2n (cf. text Eq.23) : T3

$${}_3F_2([3n-\frac{3}{2},\, n-a-\frac{1}{2},\, 2a],\, [2a+\frac{1}{2}+n,\, 3n-a-\frac{1}{2}],\, 1) \;=\;$$

$$\frac{(\dfrac{2a}{3}+\dfrac{1}{3})\,\Gamma(-3n+\dfrac{5}{2})\,\Gamma(a-n+\dfrac{1}{2})\,\Gamma(2a+\dfrac{1}{2}+n)}{\Gamma(\dfrac{1}{2}+n)\,\Gamma(\dfrac{3}{2}+a-3n)\,\Gamma(2a+2)\,\Gamma(\dfrac{3}{2}-n)}$$

excess= $2+n$

B.6  Gessel & Stanton, SIAM J. Math. Anal. 13,2,295(1982) Eq(1.9) : T7

$${}_3F_2([a,\, b,\, -n+\frac{a\,n}{b+n}],\, [1+b+\frac{a\,n}{b+n},\, a-n+1],\, 1) \;=\;$$

$$\frac{(-1)^n \sin(\pi b)\,\Gamma(1+b+\dfrac{a\,n}{b+n})\,\Gamma(-b)}{\sin(\pi a)\,\Gamma(1+a)\,\Gamma(1-n+\dfrac{a\,n}{b+n})\,\Gamma(n-a)}$$

excess= $c-(1+a+b)/2$

B.7  Prudnikov 7.4.4.20  : T1

$${}_3F_2([a,\, b,\, c],\, [2c-1,\, \frac{1}{2}a+\frac{1}{2}b+\frac{1}{2}],\, 1) \;=\;$$

$$\sqrt{\pi}\,\Gamma(c-\frac{1}{2})\,\Gamma(\frac{1}{2}a+\frac{1}{2}b+\frac{1}{2})\,\Gamma(c-\frac{1}{2}a-\frac{1}{2}b-\frac{1}{2})(\Gamma(\frac{1}{2}a)\,\Gamma(\frac{1}{2}b)\,\Gamma(c-\frac{1}{2}a)\,\Gamma(c-\frac{1}{2}b)$$
$$+\,\Gamma(\frac{1}{2}b+\frac{1}{2})\,\Gamma(c-\frac{1}{2}-\frac{1}{2}a)\,\Gamma(c-\frac{1}{2}-\frac{1}{2}b)\,\Gamma(\frac{1}{2}a+\frac{1}{2}))$$
$$\Big/(\Gamma(\frac{1}{2}b+\frac{1}{2})\Gamma(c-\frac{1}{2}-\frac{1}{2}a)\,\Gamma(c-\frac{1}{2}-\frac{1}{2}b)\,\Gamma(\frac{1}{2}a+\frac{1}{2})\,\Gamma(c-\frac{1}{2}a)\,\Gamma(c-\frac{1}{2}b)\,\Gamma(\frac{1}{2}a)\,\Gamma(\frac{1}{2}b))$$

$$\text{excess} = -\frac{1}{2}c - \frac{1}{2}a + 1 + b$$

**B.8  Prudnikov 7.4.4.20  : T2**

$$_3F_2([a, b, c], [1 + 2b, \frac{1}{2}c + \frac{1}{2}a], 1) =$$

$$\left(\left(\frac{\Gamma(\frac{1}{2}a)}{\Gamma(1 + b - \frac{1}{2}c)\Gamma(-\frac{1}{2}a + b + \frac{1}{2})\Gamma(\frac{1}{2}c)} + \frac{\Gamma(\frac{1}{2}a + \frac{1}{2})}{\Gamma(\frac{1}{2} + b - \frac{1}{2}c)\Gamma(-\frac{1}{2}a + 1 + b)\Gamma(\frac{1}{2}c + \frac{1}{2})}\right)\right.$$

$$\left.\times \sqrt{\pi}\,\Gamma(\frac{1}{2}c + \frac{1}{2}a)\Gamma(1 + 2b)\Gamma(-\frac{1}{2}c - \frac{1}{2}a + 1 + b)\right) \Big/ (2^{(-a+1+2b)}\Gamma(b+1)\Gamma(a))$$

$$\text{excess} = a - 1$$

**B.9  Prudnikov 7.4.4.20  : T7**

$$_3F_2([a, b, -b + 1], [c, 2a - c], 1) =$$

$$\left(\frac{\sqrt{\pi}\,\Gamma(a-1)\Gamma(c)\Gamma(\frac{1}{2}b + \frac{1}{2}c - \frac{1}{2})}{2^{(2a-b-c)}\Gamma(a)\Gamma(-\frac{1}{2}b + \frac{1}{2}c + \frac{1}{2})\Gamma(a - \frac{1}{2}b - \frac{1}{2}c)\Gamma(a + \frac{1}{2}b - \frac{1}{2}c - \frac{1}{2})\Gamma(c+b-1)}\right.$$

$$\left.+ \frac{\sqrt{\pi}\,\Gamma(\frac{1}{2}b + \frac{1}{2}c)\Gamma(a-1)\Gamma(c)}{2^{(2a-b-c)}\Gamma(a)\Gamma(-\frac{1}{2}b + \frac{1}{2}c)\Gamma(-\frac{1}{2}b + a - \frac{1}{2}c + \frac{1}{2})\Gamma(a + \frac{1}{2}b - \frac{1}{2}c)\Gamma(c+b-1)}\right)\Gamma(2a-c)$$

$$\text{excess} = 1 + a$$

**B.10  Prudnikov 7.4.4.20  : T8**

$$_3F_2([a, b, -b], [c, 2a - c + 1], 1) =$$

$$\frac{\sqrt{\pi}\,\Gamma(2a - c + 1)\Gamma(c)\Gamma(\frac{1}{2}b + \frac{1}{2}c)}{2^{(2a-c+1-b)}\Gamma(a - \frac{1}{2}c + 1 + \frac{1}{2}b)\Gamma(-\frac{1}{2}b + a - \frac{1}{2}c + \frac{1}{2})\Gamma(-\frac{1}{2}b + \frac{1}{2}c)\Gamma(c+b)}$$

$$+ \frac{\sqrt{\pi}\,\Gamma(2a - c + 1)\Gamma(\frac{1}{2}c + \frac{1}{2}b + \frac{1}{2})\Gamma(c)}{2^{(2a-c+1-b)}\Gamma(a - \frac{1}{2}c + \frac{1}{2} + \frac{1}{2}b)\Gamma(a - \frac{1}{2}c + 1 - \frac{1}{2}b)\Gamma(-\frac{1}{2}b + \frac{1}{2}c + \frac{1}{2})\Gamma(c+b)}$$

$$\text{excess} = -\frac{1}{2}a + 1 - \frac{1}{2}b + c$$

B.11  Prudnikov 7.4.4.22  : T1

$$_3F_2([a, b, c], [2c, 1 + \frac{1}{2}b + \frac{1}{2}a], 1) =$$

$$\left(\left(\frac{1}{2}\frac{\Gamma(c - \frac{1}{2}a)\,\Gamma(\frac{1}{2}a)}{\Gamma(c - \frac{1}{2}b)\,\Gamma(\frac{1}{2}b)} - \frac{1}{2}\frac{\Gamma(\frac{1}{2} - \frac{1}{2}a + c)\,\Gamma(\frac{1}{2}a + \frac{1}{2})}{\Gamma(\frac{1}{2} - \frac{1}{2}b + c)\,\Gamma(\frac{1}{2}b + \frac{1}{2})}\right)\Gamma(1 + \frac{1}{2}b + \frac{1}{2}a)\,\Gamma(2c)\Gamma(-\frac{1}{2}a + 1 - \frac{1}{2}b + c)\right)$$

$$\Big/ \left(\Gamma(2c - a)\,\Gamma(c)\,\Gamma(a)\,(\frac{1}{2}b - \frac{1}{2}a)\,(-\frac{1}{2}a - \frac{1}{2}b + c)\right)$$

$$\text{excess} = -\frac{1}{2}c - \frac{1}{2}b + a - \frac{1}{2}n + \frac{1}{2} + \frac{1}{2}m$$

B.12  Lewanowicz, J. Comp & Appl. Math. 86,375(1997) Generalized Watson; Eq.(2.15)  : T3

$$_3F_2([a, b, c], [\frac{1}{2}c + \frac{1}{2}b + \frac{1}{2}n + \frac{1}{2} + \frac{1}{2}m, -n + 2a], 1) =$$

$$\frac{W_{m,n}(c, -b + 2a - n, a - n, m, n)\,\Gamma(-\frac{1}{2}c - \frac{1}{2}b + a - \frac{1}{2}n + \frac{1}{2} + \frac{1}{2}m)\,\Gamma(\frac{1}{2}c + \frac{1}{2}b + \frac{1}{2}n + \frac{1}{2} + \frac{1}{2}m)}{\Gamma(-\frac{1}{2}c + \frac{1}{2}b + \frac{1}{2}n + \frac{1}{2} + \frac{1}{2}m)\,\Gamma(\frac{1}{2}c - \frac{1}{2}b + a - \frac{1}{2}n + \frac{1}{2} + \frac{1}{2}m)}$$

$$\text{excess} = b + n$$

B.13  Lewanowicz, J. Comp & Appl. Math. 86,375(1997) Generalized Watson; Eq.(2.15)  : T9

$$_3F_2([a, -a + m + 1, b], [c, -c + n + 2b + 1 + m], 1) =$$

$$\frac{W_{m,n}(a - c + n + 2b, -c - a + n + 1 + m + 2b, b, m, n)\,\Gamma(b + n)\,\Gamma(c)}{\Gamma(c - b)\,\Gamma(2b + n)}$$

$$\text{excess} = b - 2c - 2a + 3$$

B.14  Lavoie, Math. Comp., 49,179,269(1987), Eq(2)  : T2

$$_3F_2([a, b, c], [b + 2 - a, 1 + b - c], 1) =$$

$$-\frac{\Gamma(b + 2 - a)\,\Gamma(1 + b - c)\,\Gamma(\frac{1}{2}b - c - a + 2)\,\sqrt{\pi}}{\Gamma(-a + 1 + \frac{1}{2}b)\,\Gamma(\frac{1}{2}b - c + 1)\,\Gamma(\frac{1}{2}b + \frac{1}{2})\,(a - 1)\,\Gamma(-c - a + b + 2)\,2^b}$$

$$+\frac{\Gamma(b + 2 - a)\,\Gamma(1 + b - c)\,\Gamma(\frac{1}{2}b - c - a + \frac{3}{2})\,\sqrt{\pi}}{(a - 1)\,\Gamma(-c - a + b + 2)\,\Gamma(\frac{1}{2}b)\,\Gamma(\frac{1}{2}b + \frac{3}{2} - a)\,2^b\,\Gamma(\frac{1}{2} + \frac{1}{2}b - c)}$$

$$\text{excess} = a - \frac{1}{2}b - \frac{1}{2}c$$

**B.15  Lavoie, Math. Comp., 49,179,269(1987),  Eq(2)   : T3**

$$_3F_2([a,\,b,\,c],\,[1+\frac{1}{2}c+\frac{1}{2}b,\,2a-1],\,1) \;=\;$$

$$\left(\left(-\frac{1}{4}\frac{c\,(-b-1+2a)\,\Gamma(a-\frac{1}{2}b)\,\Gamma(\frac{1}{2}c+\frac{1}{2})}{\Gamma(-\frac{1}{2}c+a)\,\Gamma(\frac{1}{2}b+\frac{1}{2})} + \frac{\Gamma(1+\frac{1}{2}c)\,\Gamma(-\frac{1}{2}b+a+\frac{1}{2})}{\Gamma(-\frac{1}{2}c+a-\frac{1}{2})\,\Gamma(\frac{1}{2}b)}\right)\right.$$

$$\left.\times\Gamma(1+\frac{1}{2}c+\frac{1}{2}b)\Gamma(a-\frac{1}{2}b-\frac{1}{2}c)\,\Gamma(a-\frac{1}{2})\,2^{(c-b+2a-1)}\right) \Big/ (\sqrt{\pi}\,\Gamma(c+1)\,\Gamma(2a-b)\,(-\frac{1}{2}c+\frac{1}{2}b))$$

$$\text{excess} = b + 1$$

**B.16  Lavoie, Math. Comp., 49,179,269(1987),  Eq(2)   : T9**

$$_3F_2([a,\,-a+1,\,b],\,[c,\,-c+2+2b],\,1) \;=\;$$

$$\left(-\frac{1}{(c-1-b)\,\Gamma(\frac{1}{2}+\frac{1}{2}a-\frac{1}{2}c+b)\,\Gamma(-\frac{1}{2}c-\frac{1}{2}a+1+b)\,\Gamma(\frac{1}{2}c+\frac{1}{2}a)\,\Gamma(-\frac{1}{2}a+\frac{1}{2}c+\frac{1}{2})}\right.$$

$$\left.+\frac{1}{(c-1-b)\,\Gamma(\frac{1}{2}a-\frac{1}{2}c+1+b)\,\Gamma(-\frac{1}{2}a+\frac{1}{2}c)\,\Gamma(-\frac{1}{2}+\frac{1}{2}c+\frac{1}{2}a)\,\Gamma(\frac{3}{2}-\frac{1}{2}c-\frac{1}{2}a+b)}\right)$$

$$\times \pi\,2^{(-2b)}\,\Gamma(-c+2+2b)\,\Gamma(c)$$

$$\text{excess} = \frac{n+2-2a+ab-b}{n+1-a}$$

**B.17  Gessel & Stanton, SIAM J. Math. Anal.13,2,295(1982) Eq(1.9) : T5**

$$_3F_2([a,\,b,\,a-n-1],\,[a-\frac{n(b-1)}{a-n-1},\,a-n+1],\,1) \;=\;$$

$$\sin(\pi\,(b-a+\frac{n(-b+1)}{n+1-a}))\,(n-a)\,(bn-a+1)\,\Gamma(a+\frac{n(b-1)}{n+1-a})$$

$$\times\Gamma(2-a+\frac{(-a+1)(b-1)}{n+1-a})\,\Gamma(2-b+\frac{n(b-1)}{n+1-a}) \Big/ ((b-a+n)\,(a-1)\,\pi\,\Gamma(2+\frac{n(b-1)}{n+1-a}))$$

$$\text{excess} = -\frac{(b+n-2)(1+a-b)}{b-2}$$

**B.18  Gessel & Stanton, SIAM J. Math. Anal.13,2,295(1982) Eq(1.9) : T8**

$$_3F_2([2, a, 1-n], [b, 2 + \frac{n(-a+1)}{b-2}], 1) = \frac{(an-b-n+2)(b-1)}{(1+a-b)(b+n-2)}$$

$$\text{excess} = n + 1$$

B.19  Gessel & Stanton, SIAM J. Math. Anal. 13,2,295(1982)  Eq(1.9)
with n->-1-n    : T1

$$_3F_2([a, b, -\frac{(a-1)(b-1)}{n}], [2 + \frac{(b-1)(-a+1)}{n}, -1+a+b+n], 1) =$$

$$\frac{\sin(\pi a)(-1)^n (n-ab+a+b-1)\Gamma(n)\Gamma(-a-n+1)\Gamma(-1+a+b+n)}{\pi \Gamma(b+n)}$$

$$\text{excess} = \frac{n(n-b+a)}{a+n-1}$$

B.20  Gessel & Stanton, SIAM J. Math. Anal. 13,2,295(1982)  Eq(1.9)
with n->-1-n    : T2

$$_3F_2([a, 2, b], [1 + \frac{b(a-1)+n}{a+n-1}, 1+a+n], 1) = \frac{(n+ab-b)(a+n)}{(n-b+a)n}$$

$$\text{excess} = -\frac{3b+ab-b^2-2a-2-n+bn}{1+a-b}$$

B.21  Gessel & Stanton, SIAM J. Math. Anal. 13,2,295(1982)  Eq(1.9)
with n->-1-n    : T4

$$_3F_2([a, b, 1+a+\frac{an}{1+a-b}], [2+a, 1+a+n], 1) =$$

$$-\frac{\sin(\frac{\pi n a}{1+a-b})(1+a)\Gamma(2-b-\frac{n(-b+1)}{b-a-1})\Gamma(\frac{n(b-1)}{1+a-b})}{\sin(\pi a)\Gamma(-b+a+2+n)\Gamma(-a-n)}$$

$$\text{excess} = \frac{2a+n-2-ab+b}{a+n-1}$$

B.22  Gessel & Stanton, SIAM J. Math. Anal. 13,2,295(1982)  Eq(1.9)
with n->-1-n    : T5

$$_3F_2([a, b, a+n-1], [a-1+\frac{bn+a-1}{a+n-1}, 1+a+n], 1) =$$

$$-\sin(\pi(a+\frac{(-b+1)(a-1)}{a+n-1}))(bn+a-1)(a+n)\Gamma(1+\frac{(b-1)(-a+1)}{a+n-1})$$
$$\times \Gamma(a+\frac{n(b-1)}{a+n-1})\Gamma(2-a+\frac{(b-1)(a-1)}{a+n-1}) \Big/ ((n-b+a)\pi(a-1)\Gamma(\frac{bn+a-1}{a+n-1}+1))$$

$$\text{excess} = \frac{(n+2-b)(1+a-b)}{b-2}$$

B.23  Gessel & Stanton, SIAM J. Math. Anal. 13,2,295(1982)  Eq(1.9)
with n->-1-n    : T8

$$_3F_2([2, a, n+1], [b, 2+\frac{n(a-1)}{b-2}], 1) =$$

$$-\frac{(-1)^n \sin(\pi(b-a+\frac{n(a-1)}{b-2}))(an+b-2-n)(b-1)}{(1+a-b)(n+2-b)\sin(\frac{\pi(n+2-b)(1+a-b)}{b-2})}$$

$$\text{excess} = \frac{n(a-n-b))}{1+n-a}$$

B.24  Gessel & Stanton, SIAM J. Math. Anal. 13,2,295(1982)  Eq(1.9)
 : T2

$$_3F_2([a, 2, b], [2+\frac{(b-1)(a-1)}{a-n-1}, a-n+1], 1) = \frac{(n-a)(n-ab+b)}{(b-a+n)n}$$

$$\text{excess} = \frac{-3b-ab+b^2+2a+2-n+bn}{1+a-b}$$

B.25  Gessel & Stanton, SIAM J. Math. Anal. 13,2,295(1982)  Eq(1.9)
 : T4

$$_3F_2([a, b, 1+a-\frac{an}{1+a-b}], [2+a, a-n+1], 1) =$$

$$\frac{(1+a)\sin(\frac{\pi an}{1+a-b})\Gamma(2-b+\frac{n(b-1)}{1+a-b})\Gamma(-\frac{n(b-1)}{1+a-b})}{\sin(\pi a)\Gamma(-b+a+2-n)\Gamma(n-a)}$$

$$\text{excess} = -a+1$$

B.26   Gessel and Stanton with n->-2n (cf. text Eq.23) : T1

$$_3F_2([\frac{1}{2}+n, 3n-\frac{3}{2}, 3a+1], [2a+\frac{1}{2}+n, \frac{1}{2}+3n], 1) =$$

$$\frac{(\frac{2a}{3}+\frac{1}{3})\Gamma(-3n+\frac{5}{2})\Gamma(a-n+\frac{1}{2})\Gamma(1-a)\Gamma(\frac{1}{2}+3n)\Gamma(2a+\frac{1}{2}+n)}{\Gamma(3n-a-\frac{1}{2})\Gamma(\frac{1}{2}+n)\Gamma(\frac{3}{2}+a-3n)\Gamma(2a+2)\Gamma(\frac{3}{2}-n)}$$

$$\text{excess} = 2n$$

B.27   Gessel and Stanton with n->-2n (cf. text Eq.23) : T2

$$_3F_2([\frac{1}{2}+n, n-a-\frac{1}{2}, 2a+2-2n], [2a+\frac{1}{2}+n, \frac{3}{2}-a+n], 1) =$$

$$\frac{(\frac{2a}{3}+\frac{1}{3})\Gamma(2n)\Gamma(-3n+\frac{5}{2})\Gamma(a-n+\frac{1}{2})\Gamma(\frac{3}{2}-a+n)\Gamma(2a+\frac{1}{2}+n)}{\Gamma(3n-a-\frac{1}{2})\Gamma(\frac{1}{2}+n)\Gamma(\frac{3}{2}+a-3n)\Gamma(2a+2)\Gamma(\frac{3}{2}-n)}$$

$$\text{excess} = 3n - 3a - \frac{1}{2}$$

B.28  Gessel and Stanton with n->-2n (cf. text Eq.23) : T4

$$_3F_2([3a+1, 2a+2-2n, 2a], [2a+\frac{1}{2}+n, 2a+2], 1) =$$

$$\frac{(\frac{2a}{3}+\frac{1}{3})\Gamma(-3n+\frac{5}{2})\Gamma(a-n+\frac{1}{2})\Gamma(2a+\frac{1}{2}+n)\Gamma(3n-3a-\frac{1}{2})}{\Gamma(3n-a-\frac{1}{2})\Gamma(\frac{1}{2}+n)\Gamma(\frac{3}{2}+a-3n)\Gamma(\frac{3}{2}-n)}$$

$$\text{excess} = 2a$$

B.29  Gessel and Stanton with n->-2n (cf. text Eq.23) : T5

$$_3F_2([\frac{1}{2}+n, 3n-3a-\frac{1}{2}, 2], [\frac{1}{2}+3n, \frac{3}{2}-a+n], 1) =$$

$$\frac{(\frac{2a}{3}+\frac{1}{3})\Gamma(-3n+\frac{5}{2})\Gamma(a-n+\frac{1}{2})\Gamma(2a)\Gamma(\frac{3}{2}-a+n)\Gamma(\frac{1}{2}+3n)}{\Gamma(3n-a-\frac{1}{2})\Gamma(\frac{1}{2}+n)\Gamma(\frac{3}{2}+a-3n)\Gamma(2a+2)\Gamma(\frac{3}{2}-n)}$$

$$\text{excess} = 2a + 2 - 2n$$

B.30  Gessel and Stanton with n->-2n (cf. text Eq.23) : T6

$$_3F_2([3n-\frac{3}{2}, 3n-3a-\frac{1}{2}, 2n], [\frac{1}{2}+3n, 3n-a-\frac{1}{2}], 1) =$$

$$\frac{(\frac{2a}{3}+\frac{1}{3})\Gamma(2a+2-2n)\Gamma(-3n+\frac{5}{2})\Gamma(a-n+\frac{1}{2})\Gamma(\frac{1}{2}+3n)}{\Gamma(\frac{1}{2}+n)\Gamma(\frac{3}{2}+a-3n)\Gamma(2a+2)\Gamma(\frac{3}{2}-n)}$$

$$\text{excess} = -a + n - \frac{1}{2}$$

B.31  Gessel and Stanton with n->-2n (cf. text Eq.23) : T7

$$_3F_2([3a+1, 2, 2n], [\frac{1}{2}+3n, 2a+2], 1) =$$

$$\frac{(\frac{2a}{3}+\frac{1}{3})\Gamma(n-a-\frac{1}{2})\Gamma(-3n+\frac{5}{2})\Gamma(a-n+\frac{1}{2})\Gamma(\frac{1}{2}+3n)}{\Gamma(3n-a-\frac{1}{2})\Gamma(\frac{1}{2}+n)\Gamma(\frac{3}{2}+a-3n)\Gamma(\frac{3}{2}-n)}$$

$$\text{excess} = 3\,a + 1$$

**B.32   Gessel and Stanton with n->-2n (cf. text Eq.23) : T8**

$$_3F_2([n - a - \frac{1}{2},\, 3n - 3a - \frac{1}{2},\, 1 - a],\, [\frac{3}{2} - a + n,\, 3n - a - \frac{1}{2}],\, 1) =$$

$$\frac{(\frac{2a}{3} + \frac{1}{3})\,\Gamma(-3n + \frac{5}{2})\,\Gamma(a - n + \frac{1}{2})\,\Gamma(3a + 1)\,\Gamma(\frac{3}{2} - a + n)}{\Gamma(\frac{1}{2} + n)\,\Gamma(\frac{3}{2} + a - 3n)\,\Gamma(2a + 2)\,\Gamma(\frac{3}{2} - n)}$$

$$\text{excess} = 3\,n - \frac{3}{2}$$

**B.33   Gessel and Stanton with n->-2n (cf. text Eq.23) : T9**

$$_3F_2([2a + 2 - 2n,\, 2,\, 1 - a],\, [\frac{3}{2} - a + n,\, 2a + 2],\, 1) =$$

$$\frac{(\frac{2a}{3} + \frac{1}{3})\,\Gamma(3n - \frac{3}{2})\,\Gamma(-3n + \frac{5}{2})\,\Gamma(a - n + \frac{1}{2})\,\Gamma(\frac{3}{2} - a + n)}{\Gamma(3n - a - \frac{1}{2})\,\Gamma(\frac{1}{2} + n)\,\Gamma(\frac{3}{2} + a - 3n)\,\Gamma(\frac{3}{2} - n)}$$

$$\text{excess} = -\frac{(c - 1)(b - 1)}{a} + 1$$

**B.34   Maier[ref. 37/Chu ref. 38(see text Eq.24)**

$$_3F_2([a,\, b,\, c],\, [2 + a,\, \frac{ac + ab - a - 1 + b - bc + c}{a}],\, 1) =$$

$$\frac{(1 + a)\,\Gamma(\frac{ac + ab - a - 1 + b - bc + c}{a})\,\Gamma(\frac{-1 + b - bc + c + a}{a})}{\Gamma(\frac{ab - 1 + b - bc + c}{a})\,\Gamma(\frac{-1 + b - bc + ac + c}{a})}$$

$$\text{excess} = 2$$

**B.35   Maier[ref. 37/Chu ref. 38(cf. text Eq.24)**

$$_3F_2([a,\, b,\, c],\, [\frac{c}{2} + 1 + \frac{b}{2} + \frac{a}{2} + \frac{O}{2},\, \frac{c}{2} + 1 + \frac{b}{2} + \frac{a}{2} - \frac{O}{2}],\, 1) =$$

$$\frac{2\,(a - c - b + O)\,\Gamma(\frac{c}{2} + 1 + \frac{b}{2} + \frac{a}{2} - \frac{O}{2})\,\Gamma(\frac{c}{2} + 1 + \frac{b}{2} + \frac{a}{2} + \frac{O}{2})}{(c + a - b + O)(a - c + b + O)\,\Gamma(b + 1)\,\Gamma(c + 1)\,\Gamma(a)}$$

where

$$O = \sqrt{a^2 - 2ac - 2ab + c^2 - 2bc + b^2}$$

Appendix C

The results contained in this Appendix together with the contents of Sections 3.4 and 5 of the text are sufficient to allow the computation of any hypergeometric ($_3F_2([a, b, c], [e, f], 1)$ ) sum contiguous to any sum defined by Watson's, Whipple's or Dixon's theorems, or any of the results contiguous to those obtained from Maier's Theorems (7.1) and (7.3) [see ref. 39].

The following is the recursion in $m$ for fixed $n = N$ for the elements of the generalized Watson's theorem

$$W(m, n = N) =$$
$$-\frac{(m-1+a+b)(-3+a+b+m)(a+b-2c+3-m-2N)W(m-4, N)}{(-1+m+a+b-2c)(m-1-b+a)(-m+a-b+1)}$$
$$-((m-1+a+b)(10+2m^2+2mN+4ca+4cb-2b^2-2a^2+2Na+2Nb-8m-6N-4c)$$
$$\times W(m-2, N))/((-1+m+a+b-2c)(m-1-b+a)(-m+a-b+1))$$

(C.1)

The following is the recursion in $n$ for fixed $m = M$ for the generalized Watson's theorem.

$$W(m = M, n) = \frac{1}{2}((2c+n-1)(-3na-3nb-11n+2Mc-4ca-4cb-16c+8$$
$$+12cn+4a+4b-2M+8c^2+4n^2+2ab+Mn)W(M, n-1))$$
$$/((-1+n+c)(-n+1-2c+b)(-n+1-2c+a))$$
$$+\frac{\frac{1}{2}(2c+n-1)(2c+n-2)(a+b-2c+3-M-2n)W(M, n-2)}{(-1+n+c)(-n+1-2c+b)(-n+1-2c+a)}$$

(C.2)

This is the case $m = -1, n = 0$ for Dixon's theorem, obtained by recursion as discussed in the text.

$$X(-1, 0) = \frac{2^{(-a)}\sqrt{\pi}\,\Gamma(a-c)\,\Gamma(a-b)\,\Gamma(\frac{1}{2}a-c-b)}{\Gamma(\frac{1}{2}a+\frac{1}{2})\,\Gamma(\frac{1}{2}a-b)\,\Gamma(\frac{1}{2}a-c)\,\Gamma(a-b-c)}$$
$$+\frac{2^{(-a)}\sqrt{\pi}\,\Gamma(a-c)\,\Gamma(a-b)\,\Gamma(\frac{1}{2}a-c-b+\frac{1}{2})}{\Gamma(\frac{1}{2}a)\,\Gamma(\frac{1}{2}+\frac{1}{2}a-b)\,\Gamma(a-b-c)\,\Gamma(\frac{1}{2}+\frac{1}{2}a-c)}$$

(C.3)

This is the case $m = 2, n = 0$ for Dixon's theorem, obtained by recursion as discussed in the text.

$$X(2,0) = -2 \frac{\sqrt{\pi}\, \Gamma(\frac{7}{2} + \frac{1}{2}a - b - c)\, \Gamma(3 + a - c)\, \Gamma(3 + a - b)}{(-2+b)(-2+c)\, 2^a\, (b-1)(c-1)\, \Gamma(a+3-b-c)\, \Gamma(\frac{3}{2} + \frac{1}{2}a - c)\, \Gamma(\frac{3}{2} + \frac{1}{2}a - b)\, \Gamma(\frac{1}{2}a)}$$
$$+ ((\frac{1}{4} \frac{(2+a-2c)(-2b+a+2)(a+5-2c-2b)\sqrt{\pi}}{(-2+b)(-2+c)\, 2^a\, (b-1)(c-1)}$$
$$- \frac{1}{8} \frac{(-2b-2c+a+6)\sqrt{\pi}}{(-2+b)(-2+c)\, 2^{(-2+a)}})\Gamma(2 + \frac{1}{2}a - c - b)\, \Gamma(3 + a - b)\, \Gamma(3 + a - c))$$
$$\Big/ (\Gamma(\frac{1}{2}a + \frac{1}{2})\, \Gamma(a+3-b-c)\, \Gamma(\frac{1}{2}a + 2 - c)\, \Gamma(\frac{1}{2}a - b + 2))$$

(C.4)

Eq.(21) then gives

$$W(-1,0) = \Bigg( \frac{1}{\Gamma(\frac{1}{2} - \frac{1}{2}a + c)\, \Gamma(c - \frac{1}{2}b)\, \Gamma(\frac{1}{2}b)\, \Gamma(\frac{1}{2}a + \frac{1}{2})}$$
$$+ \frac{1}{\Gamma(\frac{1}{2} - \frac{1}{2}b + c)\, \Gamma(\frac{1}{2}b + \frac{1}{2})\, \Gamma(\frac{1}{2}a)\, \Gamma(c - \frac{1}{2}a)} \Bigg) \Gamma(-\frac{1}{2}a - \frac{1}{2}b + c)\, \Gamma(\frac{1}{2}a + \frac{1}{2}b)\Gamma(c + \frac{1}{2})\sqrt{\pi}$$

(C.5)

$$W(2,0) = \Bigg( \Bigg( - \frac{a^2 - 2ca - 2cb - 1 + 2c + b^2}{\Gamma(\frac{1}{2}a + \frac{1}{2})\, \Gamma(\frac{1}{2} - \frac{1}{2}a + c)\, \Gamma(\frac{1}{2} - \frac{1}{2}b + c)\, \Gamma(\frac{1}{2}b + \frac{1}{2})}$$
$$- \frac{8}{\Gamma(\frac{1}{2}b)\, \Gamma(\frac{1}{2}a)\, \Gamma(c - \frac{1}{2}b)\, \Gamma(c - \frac{1}{2}a)} \Bigg) \sqrt{\pi}\, \Gamma(-\frac{1}{2} - \frac{1}{2}a - \frac{1}{2}b + c)\, \Gamma(\frac{1}{2}a + \frac{1}{2}b + \frac{3}{2})$$
$$\times \Gamma(c + \frac{1}{2}) \Bigg) \Big/ ((-1 - a + b)(-a + b + 1))$$

(C.6)

The following are the forward/backward recursion formulae of the function $V_2(a, n)$ appearing in Eqs.(28) and (29) of the text.

$$V_2(a, n) = \frac{V_2(a, n-1)\left(108\,a^2 - 108\,a\,n + 54\,a + 27\,n^2 - 27\,n + 5\right)}{(12\,a + 5 - 6\,n)(12\,a + 1 - 6\,n)}$$
$$+ \frac{9\,(-3 + 2\,n)(-3\,n + 5 + 6\,a)(2\,a + 1 - n)(-3\,n + 1 + 6\,a)\,V_2(a, n-2)}{(12\,a + 5 - 6\,n)(12\,a + 11 - 6\,n)(12\,a + 1 - 6\,n)(12\,a + 7 - 6\,n)}$$

(C.7)

$$V_2(a, n-2) =$$
$$\frac{(12\,a + 5 - 6\,n)(12\,a + 11 - 6\,n)(12\,a + 1 - 6\,n)(12\,a + 7 - 6\,n)\,V_2(a, n)}{9(-3 + 2\,n)(-3\,n + 5 + 6\,a)(2\,a + 1 - n)(-3\,n + 1 + 6\,a)}$$
$$- \frac{(12\,a + 11 - 6\,n)(12\,a + 7 - 6\,n)\left(108\,a^2 - 108\,a\,n + 54\,a + 27\,n^2 - 27\,n + 5\right)V_2(a, n-1)}{9(-3 + 2\,n)(-3\,n + 5 + 6\,a)(2\,a + 1 - n)(-3\,n + 1 + 6\,a)}$$

(C.8)

The following are the recursion formula in the forward and backward directions respectively of the function $U(a, n)$ defined by the first variation of Maier's Theorem 7.3., appearing in Eq.(32) of the text.

$$U(a, n) =$$
$$- \frac{(-1 + 6\,a)(-1 + 2\,n)(4\,n - 5 - 6\,a)(18\,n^2 - 45\,n - 54\,a\,n + 42\,a^2 + 67\,a + 23)\,U(a, n-1)}{(-3\,a - 1 + n)(-1 + 2\,n - 2\,a)(4\,n - 7 - 6\,a)(6\,n - 5 - 6\,a)(6\,n - 7 - 6\,a)}$$
$$+ \frac{2\,(-1 + 2\,n)(-3 + 2\,n)(3\,n - 4 - 6\,a)(4\,n - 3 - 6\,a)(3\,n - 5 - 6\,a)(n - 2 - 2\,a)\,U(a, n-2)}{(-3\,a - 1 + n)(-3\,a - 2 + n)(6\,n - 7 - 6\,a)(4\,n - 7 - 6\,a)(-1 + 2\,n - 2\,a)(6\,n - 5 - 6\,a)}$$

(C.9)

$$U(a, n-2) =$$
$$\frac{(-3\,a - 1 + n)(-3\,a - 2 + n)(6\,n - 7 - 6\,a)(-1 + 2\,n - 2\,a)(4\,n - 7 - 6\,a)(6\,n - 5 - 6\,a)\,U(a, n)}{2(-1 + 2\,n)(-3 + 2\,n)(3\,n - 4 - 6\,a)(3\,n - 5 - 6\,a)(4\,n - 3 - 6\,a)(n - 2 - 2\,a)}$$
$$+ \frac{(-1 + 6\,a)(-3\,a - 2 + n)(4\,n - 5 - 6\,a)(18\,n^2 - 45\,n - 54\,a\,n + 42\,a^2 + 67\,a + 23)\,U(a, n-1)}{2(-3 + 2\,n)(3\,n - 4 - 6\,a)(3\,n - 5 - 6\,a)(4\,n - 3 - 6\,a)(n - 2 - 2\,a)}$$

(C.10)

This is the forward recursion for the function $V(a, n)$ defined by the second variation of Maier's Theorem 7.3., appearing in Eq.(33) of the text.

$$V(a, n) = \frac{1}{2} \frac{(3a + n - 1)(-207a - 135n + 130 + 54a^2 + 36n^2 + 108an) V(a, n-1)}{(3n - 5 + 6a)(n - 1 + 2a)(3n - 4 + 6a)}$$
$$- \frac{9}{2} \frac{(-5 + 2n)(3a + n - 1)(3a + n - 2) V(a, n-2)}{(3n - 5 + 6a)(n - 1 + 2a)(3n - 4 + 6a)}$$

(C.11)

Appendix D

The following lists 27 new evaluations of $_3F_2([a, b, c], [e, f], 1)$ obtained as discussed in the text, labelled (D.1) to (D.27) and described in the same manner as in Appendix B.

$$excess = \frac{1}{2} + a - \frac{n}{2}$$

D.1    Maier Eq. (7.1) with L=n    :T2 : T1

$$_3F_2([2a, a + \frac{1}{3} - \frac{n}{2}, a - \frac{1}{3} - \frac{n}{2}], [\frac{1}{2} + 2a - n, 3a - \frac{n}{2}], 1) =$$
$$\frac{V_2(a, n)\,\Gamma(\frac{1}{2} + a - \frac{n}{2})\,\Gamma(\frac{1}{2} + 2a - n)\,\Gamma(3a - \frac{n}{2})}{\Gamma(2a)\,\Gamma(2a + \frac{5}{6} - n)\,\Gamma(2a + \frac{1}{6} - n)}$$

$$excess = a - \frac{n}{2}$$

D.2    Maier Eq. (7.1) with L=n    :T2 : T2

$$_3F_2([2a, a - \frac{n}{2} + \frac{5}{6}, a + \frac{1}{6} - \frac{n}{2}], [\frac{1}{2} + 2a - n, 3a + \frac{1}{2} - \frac{n}{2}], 1) =$$
$$\frac{V_2(a, n)\,\Gamma(a - \frac{n}{2})\,\Gamma(\frac{1}{2} + 2a - n)\,\Gamma(3a + \frac{1}{2} - \frac{n}{2})}{\Gamma(2a)\,\Gamma(2a + \frac{5}{6} - n)\,\Gamma(2a + \frac{1}{6} - n)}$$

$$excess = 2a - \frac{1}{3}$$

D.3    Maier Eq. (7.1) with L=n    :T2 : T3

$$_3F_2([a + \frac{1}{3} - \frac{n}{2}, a - \frac{n}{2} + \frac{5}{6}, -n + \frac{1}{2}], [\frac{1}{2} + 2a - n, 2a + \frac{5}{6} - n], 1) = \frac{V_2(a, n)\,\Gamma(2a - \frac{1}{3})\,\Gamma(\frac{1}{2} + 2a - n)}{\Gamma(2a)\,\Gamma(2a + \frac{1}{6} - n)}$$

$$excess = 2a + \frac{1}{3}$$

D.4    Maier Eq. (7.1) with L=n    :T2 : T4

$$_3F_2([a - \frac{1}{3} - \frac{n}{2}, a + \frac{1}{6} - \frac{n}{2}, -n + \frac{1}{2}], [\frac{1}{2} + 2a - n, 2a + \frac{1}{6} - n], 1) = \frac{V_2(a, n)\,\Gamma(2a + \frac{1}{3})\,\Gamma(\frac{1}{2} + 2a - n)}{\Gamma(2a)\,\Gamma(2a + \frac{5}{6} - n)}$$

$$excess = -n + \frac{1}{2}$$

D.5    Maier Eq. (7.1) with L=n    :T2 : T5

$$_3F_2([2a, 2a + \frac{1}{3}, 2a - \frac{1}{3}], [3a - \frac{n}{2}, 3a + \frac{1}{2} - \frac{n}{2}], 1) = \frac{V_2(a, n)\,\Gamma(-n + \frac{1}{2})\,\Gamma(3a - \frac{n}{2})\,\Gamma(3a + \frac{1}{2} - \frac{n}{2})}{\Gamma(2a)\,\Gamma(2a + \frac{5}{6} - n)\,\Gamma(2a + \frac{1}{6} - n)}$$

$$excess = a + \frac{1}{6} - \frac{n}{2}$$

D.6    Maier Eq. (7.1) with L=n    :T2 : T6

$$_3F_2([a + \frac{1}{3} - \frac{n}{2}, 2a + \frac{1}{3}, a - \frac{n}{2}], [3a - \frac{n}{2}, 2a + \frac{5}{6} - n], 1) = \frac{V_2(a, n)\,\Gamma(a + \frac{1}{6} - \frac{n}{2})\,\Gamma(3a - \frac{n}{2})}{\Gamma(2a)\,\Gamma(2a + \frac{1}{6} - n)}$$

$$excess = a - \frac{n}{2} + \frac{5}{6}$$

D.7    Maier Eq. (7.1) with L=n    :T2 : T7

$$_3F_2([a - \frac{1}{3} - \frac{n}{2}, 2a - \frac{1}{3}, a - \frac{n}{2}], [3a - \frac{n}{2}, 2a + \frac{1}{6} - n], 1) = \frac{V_2(a, n)\,\Gamma(a - \frac{n}{2} + \frac{5}{6})\,\Gamma(3a - \frac{n}{2})}{\Gamma(2a)\,\Gamma(2a + \frac{5}{6} - n)}$$

$$excess = a - \frac{1}{3} - \frac{n}{2}$$

D.8    Maier Eq. (7.1) with L=n    :T2 : T8

$$_3F_2([a - \frac{n}{2} + \frac{5}{6}, 2a + \frac{1}{3}, \frac{1}{2} + a - \frac{n}{2}], [3a + \frac{1}{2} - \frac{n}{2}, 2a + \frac{5}{6} - n], 1) = \frac{V_2(a, n)\,\Gamma(a - \frac{1}{3} - \frac{n}{2})\,\Gamma(3a + \frac{1}{2} - \frac{n}{2})}{\Gamma(2a)\,\Gamma(2a + \frac{1}{6} - n)}$$

$$excess = a + \frac{1}{3} - \frac{n}{2}$$

D.9    Maier Eq. (7.1) with L=n    :T2 : T9

$$_3F_2([a + \frac{1}{6} - \frac{n}{2}, 2a - \frac{1}{3}, \frac{1}{2} + a - \frac{n}{2}], [3a + \frac{1}{2} - \frac{n}{2}, 2a + \frac{1}{6} - n], 1) = \frac{V_2(a, n)\,\Gamma(a + \frac{1}{3} - \frac{n}{2})\,\Gamma(3a + \frac{1}{2} - \frac{n}{2})}{\Gamma(2a)\,\Gamma(2a + \frac{5}{6} - n)}$$

$$excess = a + \frac{2}{3}$$

**D.10   Maier Thm.(7.3) with L=n variation 1   : T1**

$${}_3F_2([\frac{1}{2}, 2a+\frac{1}{3}-n, -\frac{1}{6}+n-a], [a+\frac{1}{2}, a+\frac{5}{6}], 1) = \frac{U(a,n)\,\Gamma(a+\frac{2}{3})\,\Gamma(a+\frac{1}{2})\,\Gamma(a+\frac{5}{6})}{\sqrt{\pi}\,\Gamma(3a+1-n)\,\Gamma(\frac{1}{2}+n)}$$

$$excess = a + \frac{1}{3}$$

**D.11   Maier Thm.(7.3) with L=n variation 1   : T2**

$${}_3F_2([\frac{1}{2}, 2a+\frac{2}{3}-n, \frac{1}{6}+n-a], [a+\frac{1}{2}, \frac{7}{6}+a], 1) = \frac{U(a,n)\,\Gamma(a+\frac{1}{3})\,\Gamma(a+\frac{1}{2})\,\Gamma(\frac{7}{6}+a)}{\sqrt{\pi}\,\Gamma(3a+1-n)\,\Gamma(\frac{1}{2}+n)}$$

$$excess = \frac{1}{2} + n - a$$

**D.12   Maier Thm.(7.3) with L=n variation 1   : T3**

$${}_3F_2([2a+\frac{1}{3}-n, 2a+\frac{2}{3}-n, a], [a+\frac{1}{2}, 3a+1-n], 1) = \frac{U(a,n)\,\Gamma(\frac{1}{2}+n-a)\,\Gamma(a+\frac{1}{2})}{\sqrt{\pi}\,\Gamma(\frac{1}{2}+n)}$$

$$excess = 2a + 1 - n$$

**D.13   Maier Thm.(7.3) with L=n variation 1   : T4**

$${}_3F_2([-\frac{1}{6}+n-a, \frac{1}{6}+n-a, a], [a+\frac{1}{2}, \frac{1}{2}+n], 1) = \frac{U(a,n)\,\Gamma(2a+1-n)\,\Gamma(a+\frac{1}{2})}{\sqrt{\pi}\,\Gamma(3a+1-n)}$$

$$excess = a$$

**D.14   Maier Thm.(7.3) with L=n variation 1   : T5**

$${}_3F_2([\frac{1}{2}, 2a+1-n, \frac{1}{2}+n-a], [a+\frac{5}{6}, \frac{7}{6}+a], 1) = \frac{U(a,n)\,\Gamma(a)\,\Gamma(a+\frac{5}{6})\,\Gamma(\frac{7}{6}+a)}{\sqrt{\pi}\,\Gamma(3a+1-n)\,\Gamma(\frac{1}{2}+n)}$$

$$excess = \frac{1}{6} + n - a$$

**D.15   Maier Thm.(7.3) with L=n variation 1   : T6**

$${}_3F_2([2a+\frac{1}{3}-n, 2a+1-n, a+\frac{1}{3}], [a+\frac{5}{6}, 3a+1-n], 1) = \frac{U(a,n)\,\Gamma(\frac{1}{6}+n-a)\,\Gamma(a+\frac{5}{6})}{\sqrt{\pi}\,\Gamma(\frac{1}{2}+n)}$$

$$excess = 2\,a + \frac{2}{3} - n$$

### D.16 Maier Thm.(7.3) with L=n variation 1 : T7

$$_3F_2([-\frac{1}{6}+n-a, \frac{1}{2}+n-a, a+\frac{1}{3}], [a+\frac{5}{6}, \frac{1}{2}+n], 1) = \frac{U(a,n)\,\Gamma(2\,a+\frac{2}{3}-n)\,\Gamma(a+\frac{5}{6})}{\sqrt{\pi}\,\Gamma(3\,a+1-n)}$$

$$excess = -\frac{1}{6} + n - a$$

### D.17 Maier Thm.(7.3) with L=n variation 1 : T8

$$_3F_2([2\,a+\frac{2}{3}-n, 2\,a+1-n, a+\frac{2}{3}], [\frac{7}{6}+a, 3\,a+1-n], 1) = \frac{U(a,n)\,\Gamma(-\frac{1}{6}+n-a)\,\Gamma(\frac{7}{6}+a)}{\sqrt{\pi}\,\Gamma(\frac{1}{2}+n)}$$

$$excess = 2\,a + \frac{1}{3} - n$$

### D.18 Maier Thm.(7.3) with L=n variation 1 : T9

$$_3F_2([\frac{1}{6}+n-a, \frac{1}{2}+n-a, a+\frac{2}{3}], [\frac{7}{6}+a, \frac{1}{2}+n], 1) = \frac{U(a,n)\,\Gamma(2\,a+\frac{1}{3}-n)\,\Gamma(\frac{7}{6}+a)}{\sqrt{\pi}\,\Gamma(3\,a+1-n)}$$

$$excess = a + \frac{2}{3}$$

### D.19 Maier Thm (7.3) with L=n, Variation 2 : T1

$$_3F_2([-\frac{1}{2}+n, 2\,a+n-\frac{2}{3}, -\frac{1}{6}-a], [n-\frac{1}{2}+a, -\frac{1}{6}+a+n], 1) = \frac{V(a,n)\,\Gamma(a+\frac{2}{3})\,\Gamma(n-\frac{1}{2}+a)\,\Gamma(-\frac{1}{6}+a+n)}{\Gamma(-\frac{1}{2}+n)\,\Gamma(3\,a+n)\,\sqrt{\pi}}$$

$$excess = a + \frac{1}{3}$$

### D.20 Maier Thm (7.3) with L=n, Variation 2 : T2

$$_3F_2([-\frac{1}{2}+n, 2\,a-\frac{1}{3}+n, -a+\frac{1}{6}], [n-\frac{1}{2}+a, \frac{1}{6}+a+n], 1) = \frac{V(a,n)\,\Gamma(a+\frac{1}{3})\,\Gamma(n-\frac{1}{2}+a)\,\Gamma(\frac{1}{6}+a+n)}{\Gamma(-\frac{1}{2}+n)\,\Gamma(3\,a+n)\,\sqrt{\pi}}$$

$$excess = \frac{1}{2} - a$$

### D.21 Maier Thm (7.3) with L=n, Variation 2 : T3

$$_3F_2([2a+n-\frac{2}{3}, 2a-\frac{1}{3}+n, a], [n-\frac{1}{2}+a, 3a+n], 1) = \frac{V(a,n)\,\Gamma(\frac{1}{2}-a)\,\Gamma(n-\frac{1}{2}+a)}{\Gamma(-\frac{1}{2}+n)\,\sqrt{\pi}}$$

$$excess = 2a+n$$

### D.22  Maier Thm (7.3) with L=n, Variation 2 : T4

$$_3F_2([-\frac{1}{6}-a, -a+\frac{1}{6}, a], [n-\frac{1}{2}+a, \frac{1}{2}], 1) = \frac{V(a,n)\,\Gamma(2a+n)\,\Gamma(n-\frac{1}{2}+a)}{\Gamma(-\frac{1}{2}+n)\,\Gamma(3a+n)}$$

$$excess = a$$

### D.23  Maier Thm (7.3) with L=n, Variation 2 : T5

$$_3F_2([-\frac{1}{2}+n, 2a+n, \frac{1}{2}-a], [-\frac{1}{6}+a+n, \frac{1}{6}+a+n], 1) = \frac{V(a,n)\,\Gamma(a)\,\Gamma(-\frac{1}{6}+a+n)\,\Gamma(\frac{1}{6}+a+n)}{\Gamma(-\frac{1}{2}+n)\,\Gamma(3a+n)\,\sqrt{\pi}}$$

$$excess = -a+\frac{1}{6}$$

### D.24  Maier Thm (7.3) with L=n, Variation 2 : T6

$$_3F_2([2a+n-\frac{2}{3}, 2a+n, a+\frac{1}{3}], [-\frac{1}{6}+a+n, 3a+n], 1) = \frac{V(a,n)\,\Gamma(-a+\frac{1}{6})\,\Gamma(-\frac{1}{6}+a+n)}{\Gamma(-\frac{1}{2}+n)\,\sqrt{\pi}}$$

$$excess = 2a-\frac{1}{3}+n$$

### D.25  Maier Thm (7.3) with L=n, Variation 2 : T7

$$_3F_2([-\frac{1}{6}-a, \frac{1}{2}-a, a+\frac{1}{3}], [-\frac{1}{6}+a+n, \frac{1}{2}], 1) = \frac{V(a,n)\,\Gamma(2a-\frac{1}{3}+n)\,\Gamma(-\frac{1}{6}+a+n)}{\Gamma(-\frac{1}{2}+n)\,\Gamma(3a+n)}$$

$$excess = -\frac{1}{6}-a$$

### D.26  Maier Thm (7.3) with L=n, Variation 2 : T8

$$_3F_2([2a-\frac{1}{3}+n, 2a+n, a+\frac{2}{3}], [\frac{1}{6}+a+n, 3a+n], 1) = \frac{V(a,n)\,\Gamma(-\frac{1}{6}-a)\,\Gamma(\frac{1}{6}+a+n)}{\Gamma(-\frac{1}{2}+n)\,\sqrt{\pi}}$$

$$excess = 2a+n-\frac{2}{3}$$

D.27  Maier Thm (7.3) with L=n, Variation 2 : T9

$$_3F_2([-a+\frac{1}{6}, \frac{1}{2}-a, a+\frac{2}{3}], [\frac{1}{6}+a+n, \frac{1}{2}], 1) = \frac{V(a,n)\,\Gamma(2a+n-\frac{2}{3})\,\Gamma(\frac{1}{6}+a+n)}{\Gamma(-\frac{1}{2}+n)\,\Gamma(3a+n)}$$